\documentclass{elsart}

\usepackage{amscd,amssymb}

\input xypic \xyoption{curve}
\input epsf

\newcommand{\C}[1]{{\mathcal#1}} 
\newcommand{\D}[1]{{\mathbb#1}} 

\begin{document}

\begin{frontmatter}



\title{On deformation of associative algebras and graph homology}


\author{Fusun Akman, Lucian M. Ionescu and Papa Sissokho}
\address{Department of Mathematics, Illinois State University, IL 61790-4520}

\ead{akmanf@ilstu.edu, lmiones@ilstu.edu, psissok@ilstu.edu}
\ead[url]{www.ilstu.edu/\~{ }lmiones, www.ilstu.edu/\~{}akmanf, scs.cas.ilstu.edu/\~{ }psissok}
\date{3/4/2005} 

\begin{abstract}
Deformation theory of associative algebras and in particular of Poisson algebras is reviewed.
The role of an ``almost contraction'' leading to a canonical solution
of the corresponding Maurer-Cartan equation is noted.
This role is reminiscent of the homotopical perturbation lemma,
with the infinitesimal deformation cocycle as ``initiator''.

Applied to star-products, 
we show how Moyal's formula can be obtained using such an almost contraction 
and conjecture that the ``merger operation'' provides a canonical solution 
at least in the case of linear Poisson structures.
\end{abstract}

\end{frontmatter}

\tableofcontents 



\section{Introduction}
The aim of this article is to apply perturbation techniques
to the case of the differential graded Lie algebras (DGLA) of graphs \cite{comb}
which controls the deformation theory of associative algebras (\cite{GS,DS} etc.).

Specifically, we investigate the Maurer-Cartan equation in the case of 
a differential Lie algebra in the presence of an ``almost contraction''
which leads to a ``canonical solution''.
The role of the ``merger operation'' of \cite{Kath} is unveiled,
as providing such a mapping in the well-known case of Moyal formula,
which provides a star-product in the case of a constant Poisson structure.
It is conjectured that a similar merger operation exists in the general case \ref{C:1},
where the suitable combinatorial factors are still to be determined
in a subsequent article \cite{coefs}.
The similarity with the homotopy perturbation lemma \cite{GLS} is mentioned,
to be exploited in the future work.

As a second ``improvement'' over the classical approach \cite{GS,DS},
we reduce the Maurer-Cartan equation to a Lie algebra equation,
and point out, in a special case, the role of symmetry which seems to be the key
for finding such a solution (Definition \ref{D:sym}),
a role also noted informally in the ``correction analysis'' of \cite{Kath}, p.15.

The paper is organized as follows.
We start with a brief review of Gerstenhaber theory 
of deformations of algebras \cite{GS},
phrased in the context of differential graded Lie algebras,
avoiding the Gerstenhaber pre-Lie operation.
An ``almost contraction'' (\ref{E:contr}) 
is defined and the corresponding solution is constructed.

Section \ref{S:ag} applies the above technique to the generic case of 
the DGLA of graphs.
In the constant Poisson structure case the Moyal formula is obtained in a way
which gives us hope for the general case: Conjecture \ref{C:1}.

On the other hand, since the DGLA of graphs is a 
pointed differential graded Lie algebra 
(i.e. $\partial=ad_m$, the bracket with a degree one element),
a direct proof for the associativity of the Moyal formula 
at the level of Lie algebras is provided.
It unravels a symmetry which will be studied in the general case, 
as part of the future work sketched in the concluding section.

\section{Deformation theory of associative algebras}\label{S:dt}
Given an associative algebra $(A,m)$, a {\em star product} (deformation of $m$) is 
an associative $k[[\hbar]]$-bilinear operation on $A_{\hbar}=A[[\hbar]]$ 
(\cite{DS}, p.5).
It is determined by the its values on $u,v\in A$:
\[ u\star v=m(u,v)+\hbar m_1(u,v)+\hbar^2m_2(u,v)+\cdots \]
We will recall the constraints on the coefficients imposed
by the associativity requirement.

\subsection{Maurer-Cartan equation}
Associativity of $m=m_0$ as well as of the star-product 
can be expressed conveniently using Gerstenhaber composition:
$m\circ m=0$ (\cite{Kon1}, p.9; \cite{GS}).
Let $\partial=[m, \cdot]$ be ``bracketing with $m$'', a square-zero differential,
where $[,]$ denotes the graded Lie bracket associated to the pre-Lie operation $\circ$,
where the grading is the usual shifted degree of Hochschild DGLA $g=C^\bullet(A;A)$,
so that $deg(m_i)=1$, $m_i:A\otimes A\to A$.

Grouping together coefficients of the powers of $\hbar$, 
we obtain the associativity conditions
\begin{eqnarray} && m_0\circ m_0 =0\nonumber\\
&& [m_0,m_1]=\partial m_1=0\nonumber\\
&& [m_0,m_2]+m_1\circ m_1=\partial m_2+m_1\circ m_1=0\nonumber\\
&& [m_0,m_3]+[m_1,m_2]=\partial m_3+[m_1,m_2]=0\nonumber\\
&& \vdots\nonumber\\
&& m_0\circ m_n+m_1\circ m_{n-1}+\cdots +m_{n-1}\circ m_1+m_n\circ m_0\\
&& \qquad =\partial m_n +\sum_{j,k\geq 1,\ j+k=n}m_j\circ m_k=0\nonumber\\
&& \vdots \nonumber
\end{eqnarray}
The equations are equivalent to the Maurer-Cartan equation
satisfied by the perturbation $\gamma=\star-m$ of $m$:
$$\partial\gamma+\frac12[\gamma,\gamma]=0.$$
Define trilinear maps
\[ D_n=-\sum_{j,k\ge 1,\ j+k=n}m_j\circ m_k, \quad n\ge1 \]
where the empty sum is zero.
Note that by doubling terms and using the fact $[m_j,m_k]=[m_k,m_j]$
(all $m_i$s are odd elements), we may rewrite
\[ D_n=-\frac12\sum_{j,k\ge 1,\ j+k=n}[m_j, m_k],\]
which has the advantage of involving the Lie algebra structure only, 
without making explicit use of the non-associative pre-Lie operation.
\begin{lem}
The following are equivalent:

(i) The product $\star$ is associative

(ii) $D_n=\partial m_n,\quad n\ge 1$

(iii) $[\star,\star]=0$.
\end{lem}
\begin{pf}
Regarding the equivalence between (i) and (ii),
we only need to note that 
$$[\star,\star]_{n}=\sum_{i,j\ge 0, \ i+j=n}[m_i,m_j]=2(\partial m_n-D_n).$$
\end{pf}
If the equations are satisfied up to order $r$ we say $\star$ 
is an {\em $r$-th order deformation of $m_0$}.
Then the $D_n$ satisfy the above equation up to order $r$,
i.e. $D_n$ are boundaries for $1\le n\le r$.

As a consequence the following folklore fact is obtained
(\cite{GS}; the ``simple computation'' of \cite{DS}, p.6).
\begin{lem}
Let $m_1,\dots,m_{n}$ be bilinear maps with $\partial m_1=0$.
If $D_r=\partial m_r$ are boundaries for $2\leq r\leq n$, 
then $D_{n+1}$ is a cocycle: $\partial D_{n+1}=0$. 
\end{lem}
\begin{pf}
The key point is that $\star$ is a homogeneous element of degree one (after shifting),
so that by the graded Jacobi identity
\[ [[\star,\star],\star]=0\]
the $r+1$-component vanishes
$$\sum_{i=1}^{r+1}[[\star,\star]_{i},m_{r+1-i}]=0.$$
The first $r$ terms vanish anyway,
since the assumption $D_i=\partial m_i$
is equivalent (after the ``doubling trick'') to $[\star,\star]_i=0$
(see (iii) from Lemma 1).
Therefore 
$$[[\star,\star]_{r+1},m_0]=0,$$
i.e. $[\star,\star]_{r+1}=2(\partial m_{r+1}-D_{r+1})$ is a cocycle.
Then, since $\partial m_{r+1}$ is a boundary,
$D_{r+1}$ is also a cocycle,
concluding the proof.
\end{pf}

\subsection{Obstructions}
We now review the problem of extending $r$-order deformations to $(r+1)$-order deformations
for given initial conditions: 
$$\star(0)=m, \quad \frac{d\star}{d\hbar}(0)=m_1.$$
The first extension is possible if the homology class of $D_2=-[m_1,m_1]$ is trivial.
There are no possible ``obstructions'' if $H^3(C,\partial)=Z^3/B^3$ vanishes,
where $C^m=Hom(A^m,A)$, $Z^3=ker \partial_3$ and $B_3=Im\partial_2$:
$$\diagram
0\rto & C_1 \rto^{\partial_1} & C_2\rto^{\partial_2} & C_3 \rto^{\partial_3} & ...
\enddiagram$$
On the other hand, 
the deformation is equivalent to the trivial deformation
$\star=m$ if $H^2(C,\partial)=0$.

Assume a choice of $m_2$ such that $\partial m_2=D_2$ has been made.
Then the next obstruction is the homology class of $D_3$, and so on.

Even if $H^3$ is not zero, 
an inductively defined deformation exists 
if there is an {\em almost contraction} in degree three,
i.e. a mapping $\sigma$ satisfying the equation
\begin{equation}\label{E:contr}
\sigma:\C{D}\subset Z_3\to X_2,\quad \partial\sigma+\sigma\partial=1_\C{D},
\end{equation}
where $\C{D}$ is a subspace of cocycles containing $D_n$ 
corresponding to the inductively defined $m_n$ for all $n$.

Recall that if a contracting homotopy exist globally (for $n\ge 1$):
$$\diagram
0 \rto & \dlto C_1 \dto_{Id} \rto^{\partial} & \dlto_{\sigma_2} C_2 \dto_{Id} \rto^{\partial} & 
  C_3 \dto_{Id} \rto^{\partial} \dlto_{\sigma_3} & ...\\
0 \rto & C_1 \rto^{\partial} & C_2 \rto^{\partial} & C_3 \rto^{\partial} & ... 
\enddiagram$$
then the cohomology of the complex must be trivial $H(C^\bullet,\partial)=0$.

\subsection{Almost contractions and homotopy perturbation theory}
Even if there is no contracting homotopy in degree 3,
we still have a canonical solution if there are maps $\sigma_3$ and $\sigma_4$ 
acting as an {\em almost contraction}:
$$\partial\sigma_3 D_n+\sigma_4\partial D_n=D_n,$$
which continue to satisfy this identity 
as each $D_n$ is computed out of the inductively defined $m_n$.

Indeed, if $m_n=\sigma D_n$, then 
$\partial m_n=D_n$ is equivalent to the above condition,
since $D_n$ are cocycles anyway.
In lower degrees this yields
\begin{eqnarray}
D_2&=-\frac12[m_1,m_1],\\
m_2&=\sigma D_2=-\frac12\sigma([m_1,m_1]),\\
D_3&=-\frac12([m_2,m_1]+[m_1,m_2])=\frac12[m_1,\sigma[m_1,m_1]],\\
m_3&=\sigma D_3=\frac12\sigma[m_1,\sigma[m_1,m_1]].
\end{eqnarray}
Define $t=\sigma\circ ad_{m_1}$ and $\hat{m}_1^{n+1}=t^n(m_1), n\ge 0$.
Then we have
\begin{eqnarray}
D_4&=-\frac12([m_3,m_1]+[m_2,m_2]+[m_1,m_3]),\\
m_4&=\sigma D_4=\hat{m}_1^4-\frac12\sigma([\hat{m}_1^2,\hat{m}_1^2]).
\end{eqnarray}
It is natural to investigate the conditions under which such a 
``minimal procedure'' with ``initiator'' $t$ and cocycle $m_1$ exists.
Its interpretation from the perspective of the Homotopical Perturbation Lemma
(\cite{HS}, p.10) will be considered elsewhere. 

A case when such a procedure is successful 
is the one of the Moyal star-product
$$\star=exp(\hbar m_1),$$
as it will be explained next,
at the level of graphs.

\section{Application to graphs}\label{S:ag}
Let $\C{G}_{n,m}$ be the set of {\em orientation classes} 
of {\em Lie admissible edge labeled graphs} of \cite{MP}, p.3,
corresponding to {\bf linear Poisson structures} (see also \cite{comb}).
An element $\Gamma\in \C{G}_{n,m}$ is a directed graph with $n$ internal vertices, 
$m$ labeled boundary vertices $1,2...,m$ (left to right in figures), 
such that each internal vertex is trivalent with exactly two descendants.
The corresponding two outgoing arrows will be labeled left/right,
defining the {\em orientation class} of the graph $\Gamma$ 
up to a ``negation'' of the edge labeling in any two internal vertices \cite{MP}.
The orientation class of graph embedded in the plane
will be determined by the positive orientation of the plane.
The corresponding (graded) space is denoted by $\C{G}=\cup \C{G}_m$,
where $\C{G}_m=\cup_{n\in \D{N}} \C{G}_{n,m}$.
Let $C=k\C{G}$ be the quotient of the DGLA of graphs of \cite{comb},
with pre-Lie composition $\circ$ and 
differential $\partial=[b_0, \cdot]$, where $b_0\in \C{G}_{0,2}$,
by the ideal generated by the Jacobi identity (\ref{E:J}).

The initial conditions of the ``universal'' deformation problem 
are $m_0=b_0$ and $m_1=b_1$, where

$$
\begin{xy}
,(-5,3)*{b_0=\ };
,(5,-4)*{\circ};
,(5,7);(15,-4)*{\circ}
\end{xy}
\hspace{.5in}
\begin{xy}
,(-5,3)*{b_1=\ }
,(10,7)*{\bullet};(5,-4)*{\circ}**\dir{-} ?>*\dir{>}
,(10,7);(15,-4)*{\circ}**\dir{-}?>*\dir{>}
\end{xy}
\qquad.$$

Recall that $b_0\circ b_0=0$ and $[b_0,b_1]=0$ (\cite{comb} p.13).

The first possible obstruction is the homology class of 
$$D_2=-b_1\circ b_1=-(t_2^R-t_2^L+c_2^L-c_2^R)$$
where
$$
\ \ c_2^R= \ \
\begin{xy}
,(5,7)*{\bullet};(-5,-4)*{\circ}**\dir{-} ?>*\dir{>}
,(5,7);(15,-4)*{\circ}**\dir{-}?>*\dir{>}
,(15,7);(5,-4)*{\circ}**\dir{-} ?>*\dir{>}
,(15,7)*{\bullet};(15,-4)*{\circ}**\dir{-} ?>*\dir{>}
\end{xy}
\ \  \mbox{ and }\ c_2^L= \ \ 
\begin{xy}
,(5,7)*{\bullet};(-5,-4)*{\circ}**\dir{-} ?>*\dir{>}
,(-5,7);(5,-4)*{\circ}**\dir{-}?>*\dir{>}
,(5,7);(15,-4)*{\circ}**\dir{-} ?>*\dir{>}
,(-5,7)*{\bullet};(-5,-4)*{\circ}**\dir{-} ?>*\dir{>} 
\end{xy}
$$
and the graphs $t_2^R, t_2^L$ are depicted in the LHS of the 
following diagram representing the Jacobi identity $t_2^R-t_2^L=c_2$
\begin{equation}\label{E:J}
\begin{xy}
,(0,2)*{\bullet};(-5,-4)*{\circ}**\dir{-} ?>*\dir{>}
,(0,2);(5,-4)*{\circ}**\dir{-}?>*\dir{>}
,(4,7);(0,2)**\dir{-} ?>*\dir{>}
,(4,7)*{\bullet};(15,-4)*{\circ}**\dir{-} ?>*\dir{>}
\end{xy}
\ \  -\ \ 
\begin{xy}
,(4,7)*{\bullet};(-5,-4)*{\circ}**\dir{-} ?>*\dir{>}
,(9,2);(5,-4)*{\circ}**\dir{-}?>*\dir{>}
,(4,7);(9,2)**\dir{-} ?>*\dir{>}
,(9,2)*{\bullet};(15,-4)*{\circ}**\dir{-} ?>*\dir{>} 
\end{xy}
\ \ =\ \ 
\begin{xy}
,(-1,2)*{\bullet};(-5,-4)*{\circ}**\dir{-}?>*\dir{>}
,(4,7)*{\bullet};(5,-4)*{\circ}**\dir{-}?>*\dir{>}
,(4,7);(-1,2)**\dir{-} ?>*\dir{>}
,(-1,2);(15,-4)*{\circ}**\crv{(0,4)} ?>*\dir{>}
\end{xy}
\end{equation}
Using this identity, $D_2$ simplifies to $D_2=c_2^R-c_2^L-c_2$
(for additional details, see \cite{comb} p.16; \cite{CI}, p.20).

\subsection{Candidates for almost contractions}
We claim that an almost contraction as needed earlier is the ``merger operation''
(\cite{comb}, p. 10; see also \cite{Kath}, p.17):
\begin{eqnarray}
\sigma_i(\Gamma)&=\Gamma/(i(i+1)),\quad \Gamma\in\C{G}_{n,m}\label{E:merger}\\
\sigma(\Gamma)&=\frac1{2(2^n-2)} \sum_{i=1}^{m-1} (-1)^{i-1}\sigma_i(\Gamma), 
\end{eqnarray}
where the quotient graph from the RHS of (\ref{E:merger}) is obtained by merging
the $i$th and $i+1$st boundary points.
If a non-admissible graph emerges after the merger, 
the result is considered to be zero.

For example, we have $\sigma(c_2^R)=\sigma_1(c_2^R)=\frac14 b_1^2$ (similarly $\sigma(c_2^L)=-\frac12 b_1^2$):
$$
\begin{xy}
,(5,7)*{\bullet};(-5,-4)*{\circ}**\dir{-} ?>*\dir{>}
,(5,7);(15,-4)*{\circ}**\dir{-}?>*\dir{>}
,(15,7);(5,-4)*{\circ}**\dir{-} ?>*\dir{>}
,(15,7)*{\bullet};(15,-4)*{\circ}**\dir{-} ?>*\dir{>}
\end{xy}
\quad
\begin{xy}
,(-5,3)*{\mapsto}
,(5,7)*{\bullet};(5,-4)*{\circ}**\dir{-} ?>*\dir{>}
,(5,7);(15,-4)*{\circ}**\dir{-}?>*\dir{>}
,(15,7);(5,-4)**\dir{-} ?>*\dir{>}
,(15,7)*{\bullet};(15,-4)*{\circ}**\dir{-} ?>*\dir{>}
\end{xy}
$$
We will investigate the above claims in the special cases of 
constant and linear Poisson structures.

\subsection{Constant Poisson structures}
As an example we derive Moyal's formula
along the previous lines using the ``merger of legs''
as an almost contracting operation.

The benefit of having a Poisson structure with constant coefficients is
that a graph with an arrow landing on an internal vertex evaluates to zero under
Kontsevich rule $B(\Gamma)=\C{U}_\Gamma(\alpha^{\wedge n})$
where $\Gamma\in\C{G}_{n,m}$ (\cite{Kon1}, p.23, p.28):

Therefore
$$\begin{xy}
,(-10,4)*{\Gamma=b_1^n =};
,(6,8)*{\overbrace{
\begin{xy}
,(8,7)*{\bullet};
,(10,7)*{\bullet};
,(12,7)*{\bullet};
,(5,7)*{\bullet};(5,-4)*{\circ}**\dir{-} ?>*\dir{>}
,(5,7);(15,-4)*{\circ}**\dir{-}?>*\dir{>}
,(15,7);(5,-4)*{\circ}**\dir{-} ?>*\dir{>}
,(15,7)*{\bullet};(15,-4)*{\circ}**\dir{-} ?>*\dir{>}
\end{xy}
}^{\mbox{ $n$ wedges }}};
\end{xy}
$$
is the unique graph in $\C{G}_{n,2}$ not in the kernel of $B$.

In particular, the Jacobi identity (\ref{E:J}) is automatically satisfied,
since all the terms evaluate to zero under Kontsevich rule
$$B(t_2^R)=B(t_2^L)=B(c_2)=0.$$
\begin{lem}
For any $i,j\ge 0$ we have
$$\sigma([b_1^i,b_1^j])=-\frac1{2^{i+j-1}-1}b_1^{i+j},$$
where $b_1^n\in\C{G}_{n,2}, n\ge1$, with the natural orientation.
\end{lem}
\begin{pf}
It is enough to note that the only term of $b_1^i\circ_1 b_1^j$ not vanishing
after the application of $\sigma$,
is the one for which all $i$ of the left legs of $b_1^i$
land on the left boundary point of $b_1^j$, since otherwise all consecutive boundary points
are ``bridged'' by some $b_1$,
and therefore the term vanishes under the merger operation
$$\sigma(b_1^i\circ_1 b_1^j)=-\frac1{2(2^{i+j}-2)}b_1^{i+j}.$$
\end{pf}
It follows that $m_2=\sigma D_2=b_1^2/2$ and in general, we have
\begin{lem}
If $m_0=b_0, m_1=b_1$ and $m_n=\sigma D_n, n\ge 2$ then $\forall n, m_n=b_1^n/n!$.
\end{lem}
\begin{pf}
Assuming inductively that $m_k=b_1^k/k!$ for $1\le k\le n-1$, then
\begin{eqnarray}
m_n=\sigma D_n&=-\frac12\sum_{i+j=n, \ i,j\ge1}\sigma\left(\left[\frac{b_1^i}{i!},\frac{b_1^j}{j!}\right]\right)\\
&=\left(-\frac12\right)\left(-\frac1{2^{n-1}-1}\right)b_1^n\sum_{i+j=n, \ i,j\ge1}\frac1{i!} \frac1{j!}=\frac{b_1^n}{n!}.
\end{eqnarray}
\end{pf}
Now since the Moyal formula provides an associative product
$$*=e^{b_1 h},\quad [*,*]=0,$$
$D_n=\partial m_n$ are boundaries and therefore,
together with $m_n=\sigma D_n$, it implies that $\sigma$ is an almost contraction
for the inductively defined $m_n=\sigma D_n$, starting with the cocycle $m_1$:
$$\partial\sigma D_n+\sigma\partial D_n=D_n, \quad n\ge 2.$$
This, of course, amounts to $\partial\sigma D_n=D_n$,
which in turn is equivalent to the original equation in degree $n$.
Therefore we will give a direct proof that the above star-product is associative,
in order to better understand the combinatorics involved.
In contrast with the previous more general approach,
we will take advantage of the fact that the differential $\partial$ is defined
as a Lie bracket,
and focus on the Lie algebra structure.
\begin{prop}\label{P:star0}
$$[*,*]=0.$$
\end{prop}
\begin{pf}
The $n-th$ homogeneous degree of the above equation is:
\begin{equation}\label{eq:nth}
\sum_{i+j=n, \ i,j\ge 0}[m_i,m_j]=0,\quad m_k=b_1^k/k!.
\end {equation}
To prove it we will start by determining the structure coefficients
of the Lie bracket.
In order to isolate the combinatorial factors 
from the Lie algebra structure constants,
it is better to adopt a basis with elements of the form
$\Gamma/|Aut(\Gamma)|$.

Consider $\{B_n=b_1^n/n!\}$ as a basis in $k\C{G}_{\bullet,2}$.
Incidentally, the solution of $[Z,Z]=0$ is therefore the corresponding ``integral''
$*=\sum_n B_n$.

Consider the graphs $\Gamma_1,\Gamma_2,\Gamma_3\in \C{G}_{1,3}$, defined as follows:
$$
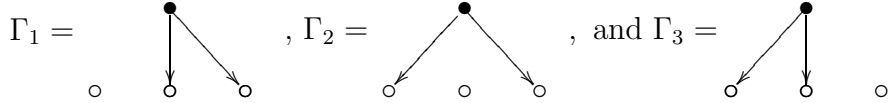
\begin{xy}
,(-10,4)*{\Gamma_1=\quad };
,(5,7)*{\bullet};(5,-4)*{\circ}**\dir{-} ?>*\dir{>}
,(5,7);(15,-4)*{\circ}**\dir{-}?>*\dir{>}
,(15,7);(5,-4)*{\circ}
,(15,-4)*{\circ}
,(-5,-4)*{\circ};
\end{xy}
\ \   
\begin{xy}
,(-12,4)*{\mbox{ , }\Gamma_2=\quad };
,(5,7)*{\bullet};(-5,-4)*{\circ}**\dir{-} ?>*\dir{>}
,(5,7);(15,-4)*{\circ}**\dir{-}?>*\dir{>}
,(5,-4)*{\circ};
\end{xy}
\ \   
\begin{xy}
,(-6,4)*{,\mbox{ and }\Gamma_3=\ }; (5,-4)*{\circ}
,(5,7);(15,-4)*{\circ}
,(15,7);(5,-4)*{\circ}**\dir{-} ?>*\dir{>}
,(15,7)*{\bullet};(15,-4)*{\circ}**\dir{-} ?>*\dir{>}
,(25,-4)*{\circ};
\end{xy}
$$
Then 
$$\{\Gamma_{rst}=(\Gamma_1^r/r!)(\Gamma_2^s/s!)(\Gamma_3^t/t!)\}_{r,s,t\ge 0}$$
is a basis in $k\C{G}_{\bullet,3}$ and 
$$\forall i,j\ge0\quad [B_i,B_j]=\sum_{r+s+t=i+j} C_{(i,j)}^{(r,s,t)} \Gamma_{rst}.$$
To compute the coefficients $C_I^J$ of $\Gamma_J$, where $I=(i,j)$ and $J=(r,s,t)$,
consider $b_1^i\circ_1 b_1^j$ first
and note that when splitting the $i$-left legs of $b_1^i$ to make them land
on the first two boundary points of $b_1^j$, 
the only graphs $\gamma=\Gamma_1^r\Gamma_2^s\Gamma_3^t$ that are involved 
are those for which $r+s=i, t=j$.
\begin{enumerate}
\item If $r+s=i$ and $t=j$ then $b_1^i\circ_1b_1^j$ contributes $i!/(r!s!)$ to $\gamma$, 
thus $C_I^J=1$.

\item If $r=i$ and $s+t=j$ then $b_1^j\circ_2b_1^i$ contributes $-j!/(s!t!)$ to $\gamma$, 
thus $C_I^J=-1$. 

\item If $r+s=j$ and $t=i$ then $b_1^j\circ_1b_1^i$ contributes $j!(r!s!)$ to $\gamma$,
thus $C_I^J=1$.

\item If $r=j$ and $s+t=i$ then $b_1^i\circ_2b_1^j$ contributes $-i!/(s!t!)$ to $\gamma$,
thus $C_I^J=-1$.

\item  If none of the above cases hold then $\gamma$ is not present in $[B_i,B_j]$, thus $C_I^J=0$.
\end{enumerate}

In conclusion we have the following lemma.

\begin{lem}\label{L:B}
\begin{eqnarray}\label{eq:[bi,bj]}
\forall i,j\geq0, \quad [B_i,B_j]
&=&\sum_{r+s=i,\ t=j}\Gamma_{(r,s,t)}
-\sum_{r=i,\ s+t=j}\Gamma_{(r,s,t)}\cr
&\ & +\sum_{r+s=j,\ t=i}\Gamma_{(r,s,t)}
-\sum_{r=j,\ s+t=i}\Gamma_{(r,s,t)}.
\end{eqnarray}
\end{lem}
To understand the algebraic reason for the cancelation better,
define the following codifferential (dual to addition in some sense):
$$\delta(i,j)=\sum_{r+s=i,\ t=j}(r,s,t)-\sum_{r=i,\ s+t=j}(r,s,t).$$
Then bracket in Lemma \ref{L:B} is its symmetrization:
$$[B_i,B_j]=<\Gamma,\delta(i,j)+\delta(j,i)>, \quad 
\Gamma(r,s,t)=W_{i,j}^{(r,s,t)}\Gamma_{(r,s,t)},$$
where $\Gamma$ is the linear operator extending the function 
defined on the corresponding domain in the $(r,s,t)$-space.
The $W(r,s,t)=1$ are the ``true coefficients'' of the Lie bracket,
without the grading sign built into $\circ$,
which is independent of the particular case under consideration.

For a geometric viewpoint of the ``integration domain'',
consider the 3-simplex $0\le r,s,t\le n$, where $n=i+j$ is fixed.
Then $\{(r,s,t)|r+s+t=i+j=n\}$ is the front face,
$r+s=i, t=j$ defines a segment parallel to the $rs$-plane and
$r=i, s+t=j$ defines a segment parallel to the $st$-plane, 
both contained in the front face and having $(i,0,j)$ as common point.

When summing over $(i,j)$, $i+j=n$, both segments swipe the front face
$$\{r+s=i, t=j, i+j=n\}=\{r+s+t=n\}=\{r=i, s+t=j, i+j=n\}.$$
Now, due to the opposite signs, there is an overall cancelation:
\begin{lem}
$$\sum_{i+j=n, \ i,j\ge 0}\delta(i,j)=0.$$
\end{lem}
As a corollary, (\ref{eq:nth}) holds true,
concluding the proof of the Proposition.
\end{pf}
Note that the proof of the proposition does not depend on the values $W(r,s,t)$, 
but rather on a certain symmetry of the basis elements involved in the Lie bracket.
\begin{defn}\label{D:sym}
The {\em antipodal map} of the DGLA of graphs is (\cite{comb}):
$$S(\Gamma)=(-1)^m \Gamma^t,\quad \Gamma\in \C{G}_{n,m}$$
where $\Gamma^t$ is the {\em transpozed graph}, 
i.e. the graph obtained by reversing the order on the boundary points.
\end{defn}
For example $S(b_1)=-b_1^t=b_1$, since they define the same orientation class.
\begin{lem}
The antipodal map is a pre-Lie morphism:
$$S(\Gamma_1\circ\Gamma_2)=S(\Gamma_1)\circ S(\Gamma_2),$$
and therefore an involution of the Lie algebra of graphs.
\end{lem}
The role of the symmetrization of a star-product was already noted 
in \cite{Kath} and \cite{comb}.
\begin{rem}
If we define:
$$\delta(n)=\sum_{i+j=n, \ i,j\ge 0}(i,j)$$
then the previous lemma says that $\delta^2=0$,
i.e. $\delta$ is indeed a codifferential.

Note also that $\delta$ is associated with the asymmetric operation:
$$\{\Gamma_1,\Gamma_2\}=\Gamma_1\circ_1\Gamma_2-\Gamma_2\circ_2\Gamma_1, 
\quad \Gamma_i\in G_{\bullet,2}.$$
Its properties will be investigated elsewhere.
\end{rem}
As a second example we will consider the case of linear Poisson structure.

\subsection{Linear Poisson structures}
Explicit star-products for linear Poisson structures (e.g. dual of a Lie algebra)
were known to exist since \cite{Gutt,MP,Gutt2}.

In this case the graphs not in the kernel of the Kontsevich rule
are products of tree-like graphs, since at most one arrow may land on
internal vertex in order to have a non-zero contribution.

A candidate for an almost contraction is 
the ``merger operation'' (\ref{E:merger}). 
\begin{lem}
$\sigma$ is a homological differential.
$$\sigma^2=0.$$
\end{lem}
\begin{pf}
Indeed, if $j\ge i$ then $\sigma_j\circ\sigma_i=\sigma_i\circ\sigma_{j+1}$
and the opposite sign of the two terms yields a pairwise cancelation
as usual.
\end{pf}
Specializing to degrees two and three we obtain
$$\sigma_2(\Gamma)=-1/(2^{n-1}-1)\Gamma/b_0, \ \Gamma\in \C{G}_{n,2},$$
$$\sigma_3(\theta)=-1/(2^{n-1}-1)(\theta/(12)-\theta/(23)),\ \theta\in \C{G}_{n,3}.$$
At present the relation between the two differentials $\sigma$ and $\partial$
is not clear.

Some elementary facts are recorded next.
\begin{lem}
For any graph $\Gamma\in \C{G}_{n,1}$ we have
$$(\partial\Gamma)/b_0=2^{i-1}\Gamma,$$
where $i$ is the number of edges landing on the unique boundary vertex.
\end{lem}
For Bernoulli graphs $b_n$ \cite{comb}, p.5, we have the following.
\begin{lem}
(i) $\partial\sigma_2(b_n^L)=0$
(ii) $\sigma_3\partial(b_n^L)=2^{n-1}b_n^L-S_R(b_n^L)/b_0^R$
where $S_R$ (respectively $S_L$) splits in all non-trivial ways
the arrows landing on $L$ ($R$).
\end{lem}
\begin{conj}\label{C:1}
A canonical solution is defined inductively by $Z_n=\sigma D_n$.
\end{conj}
Although stated in the context of linear Poisson structures,
we believe that the above conjecture holds in general,
with the appropriate combinatorial coefficients for the merger
operations $\sigma_i$, to be discussed elsewhere \cite{coefs}.

\section{Conclusions}\label{S:c}
We showed that Maurer-Cartan equation can be solved provided that there is 
an almost contraction.
This is reminiscent of the homotopy perturbation lemma 
with the infinitesimal cocycle as ``initiator'' \cite{GLS,HS}.
As an application to star-products,
the Moyal's formula was obtained in this way.

It is conjectured that the ``merger operation'',
which is a homology differential,
provides such an almost contraction at least in the case of linear Poisson structures,
leading to a canonical star-product.
Further investigations will be reported in a forthcoming article \cite{coefs}.




\end{document}